\documentclass[a4in]{article}
\usepackage{a4wide, latexsym, graphics, amssymb}
\newtheorem{corollary}{Corollary}

\newtheorem{theorem}[corollary]{Theorem}

\newtheorem{rem}[corollary]{Remark}
\newcommand{\Prob} {{\bf P}}
\newcommand{\R}{\mathbb{R}}
\newcommand{\C}{\mathbb{C}}

\newcommand{\expect}{{\bf E}}
\def\buildrel#1\over#2{\mathrel{\mathop{\kern
0 pt #2}\limits^{#1}}}
\def\build#1_#2^#3{\mathrel{\mathop{\kern 0
pt#1}\limits_{#2}^{#3}}}
\title {
Critical exponents for two-dimensional percolation
}
\author {Stanislav Smirnov\footnote {KTH Stockholm and KVA}
\and 
Wendelin Werner\footnote{Universit\'e Paris-Sud and IUF}
}
\date {December 2001}
\begin {document}

\maketitle

\begin {abstract}
We show how to combine Kesten's scaling relations,
the determination of critical exponents associated to the
stochastic Loewner evolution process by Lawler, Schramm, 
and Werner, and Smirnov's proof of Cardy's formula,
in order to determine the existence and value of critical exponents
associated to percolation on the triangular lattice.
\end {abstract}

{\bf MSC class.:} 60K35, 82B27, 82B43

\section {Introduction}

The goal of the present note is to review and clarify
the consequences
of recent papers and preprints 
concerning the existence and values of critical exponents for
site percolation on the triangular lattice.

Suppose that $p \in (0,1)$ is fixed. Each vertex of the
triangular lattice (or equivalently each hexagon in the
honeycomb lattice) is open (or colored blue) with 
probability $p$ and 
closed
(or colored yellow) with
probability $1-p$, independently of each other. It is now well-known
(and due to Kesten and
Wierman, 
see the textbooks \cite {Kbook,G}) that when $p \le 1/2$, there is
almost surely no infinite cluster of open vertices, while if $p>1/2$,
there is a.s. a unique such infinite cluster and the probability
$\theta (p)$ that the origin
belongs  to this infinite cluster is then positive.
Arguments from theoretical physics predicted that
when $p$ approaches the critical value $1/2$ from
above, $\theta (p) $ behaves roughly like
$(p-1/2)^{5/36}$.
This number $5/36$
is one of several critical exponents that are supposed
to be independent of the considered planar lattice and that are
describing the behaviour of percolation near its
critical point $p=p_c$ ($p_c = 1/2$ for this particular model).
Our goal in the present paper is to point out that this
result, as well as other related statements,
is a consequence of the combination of
various papers:

\begin {itemize}
\item
In \cite {Kpaper}, Kesten has shown that in order to
understand the behaviour of percolation near its
critical
point (and in particular existence and values of certain
critical exponents), it is sufficient to study
 what happens
at the critical point i.e. here when $p=p_c=1/2$.
In particular, many results would follow from 
the existence and the values  of
the exponents  describing the decay when $R$ goes to infinity
of the probabilities (at $p=1/2$)
 of the events $A_R^1$
and $A_R^2$ that if we restrict the percolation
to the disc of radius $R$,
there exist one (respectively two disjoint)
blue clusters
joining the vicinity of the origin (say that are at distance less
than two of the origin)
to the
circle of radius $R$.
\item
In \cite {Sch}, Schramm defined a family of
random evolution processes based on Loewner's equation,
and pointed out that one of them (the stochastic
Loewner evolution process with
parameter $6$, referred to as $SLE_6$ in the
sequel) is the only possible conformally
invariant scaling limit of discrete critical
percolation cluster interfaces.
In a series of papers \cite {LSW1, LSW2, LSW3, LSWprep}, Lawler, Schramm,
and Werner have  derived various  properties of
$SLE_6$ (relation between radial and chordal processes,
locality property etc),
  computed
critical exponents associated to $SLE_6$, and then used these
exponents to determine the value of Brownian exponents 
(for instance, the Hausdorff dimension of the planar Brownian 
frontier is a.s. $4/3$).

\item
In \cite {S, Sprep}, Smirnov proved that indeed, critical
site  percolation on the triangular lattice has a conformally
invariant scaling limit when the mesh of the lattice
tends to zero, and in particular, that the
discrete cluster interfaces converge to this stochastic
Loewner evolution process.

\end {itemize}

In the present paper, 
we will outline how one can combine all these
results to show that the critical exponents for discrete
site percolation on the triangular lattice are those
predicted in the physics literature.
All the results that we prove in the present paper 
have been conjectured (on the basis
of numerics and heuristics) 
and predicted  
(using Coulomb Gas methods, Conformal Field Theory, or Quantum Gravity) 
by physicists. 
See e.g., Den Nijs \cite {DN}, Nienhuis et al. \cite {N, N2},
Pearson \cite{P}, 
Cardy \cite {Ca}, Sapoval-Rosso-Gouyet \cite {SRG},
Grossman-Aharony \cite {GA}, Duplantier-Saleur \cite {DS},
 Cardy \cite {Cardy},
Duplantier \cite {D},
Aizenman-Duplantier-Aharony \cite {ADA} and the references in these
papers.
There exists a vast theoretical physics literature on this subject and
we do not claim that this list covers all important  
contributions to it.

Note the description of the scaling limit via $SLE_6$ enables 
also to derive some results that have not appeared in the 
physics literature. For instance,  
an analogue of Cardy's formula ``in the bulk'' \cite {Sch2}, 
or a description of the so-called backbone exponent \cite {LSWprep}.

\section {Kesten's scaling relations}

Let 
us now introduce some 
notation.
Denote by $N$ the cardinality of 
the cluster $C$ containing the origin. 
Recall that 
$$ \theta (p) = \Prob_p [ N = \infty ]
$$
where 
 $\Prob_p, \expect_p$ corresponds to site percolation on the 
triangular lattice with 
parameter $p$. Let 
$$ \chi (p) = \expect_p [ N \, 1_{N<\infty} ].$$
 This corresponds to the average cardinality of 
finite clusters.
Let 
$$
\xi (p) =
\left[ 
\frac { \expect_p [ \sum_{y \in C} |y|^2 1_{N<\infty} ] 
}{\chi (p)} 
\right]^{1/2} .$$
This is the so-called correlation length corresponding to 
the ``typical 
radius'' of a finite cluster. Other
 definitions of correlation length are 
possible, see e.g., \cite {G,Kpaper,CCF}; for 
instance, define 
 $\xi^* (p)$
by the relation 
$$
\Prob [ 0 \hbox { is connected to } x \hbox { by a finite cluster}] 
= \exp \{ - x/ \xi^* (p)  + o(x) \}
\hbox { when } x \to +\infty.$$ 
If $x$ tends to infinity along some fixed direction, 
existence of $\xi^*$ (dependent on the
chosen direction) easily follows. 
So to make the definition rigorous, one can assume for instance 
that $x\in\R$.
On the other hand, the proof  shows that the 
asymptotic behaviour of $\xi^*(p)$ in the neighbourhood of $p_c=1/2$
is independent of the chosen direction (i.e., (iv) holds for 
all given directions).

\begin {theorem}[Behaviour near the critical point]
\label {main}
~\\
\begin {itemize}
\item {(i)}
When $p \to 1/2+ $, 
$$\theta(p) = (p-1/2)^{5/ 36 + o(1)}.$$
\item {(ii)}
When $p \to 1/2$
$$
\chi (p) = (p-1/2)^{-43/18 + o(1)}.$$
\item {(iii)}
When $p \to 1/2$,
$$
\xi (p) = ( p-1/2)^{-4/3 + o(1)}.$$
\item {(iv)}
When $p \to 1/2$,
$$
\xi^* (p) = (p-1/2)^{-4/3 + o(1)}.$$
\end {itemize}
\end {theorem}

It has been shown by Kesten in \cite {Kpaper} (using 
also 
the remark following Lemma 8 in \cite {Kpaper}; 
it can be shown 
that (iv) holds once (iii) holds, 
using the estimates in \cite {Kpaper}, see
 also section 3 in \cite {CCF})
that all these results hold provided that: 
When $p=1/2$ 
and $R \to \infty$,
\begin {equation}
\label{onearm}
\Prob [ A_R^1 ] = R^{-5/48 + o(1)}
\end {equation}
and
\begin {equation}
\label {4arms}
\Prob [ A_R^2 ] = R^{-5/4 + o(1)}.
\end {equation}
The relations between these two 
critical exponents and those 
appearing in Theorem \ref {main} are sometimes
known as scaling relations.

Relation (\ref {onearm}) is proved in \cite {LSWprep}, combining the 
computation of an exponent for $SLE_6$ with the results 
of \cite {S, Sprep}.
In the rest of this paper, we fix $p=1/2$ and we shall see 
how (\ref {4arms})
 and other closely related results (we will also briefly 
discuss (\ref {onearm})) follow from 
the combination of 
\cite {S, Sprep} with \cite {LSW1,LSW2,LSW3,LSWprep}.
Independently, Yu Zhang \cite {Z} has recently anounced a proof of 
(\ref {4arms}) and of (iii) in the above
Theorem, probably using similar thoughts and arguments as
those that we shall present here.

\section {Half-plane exponents}

In this section, we are going to study the decay when $R \to \infty$,
of the probability of the events 
that there exist $j$ disjoint blue paths
that stay in the upper half-plane, start at the vicinity of
the origin, and reach distance $R$.
Besides being of independent interest,
this section serves as a model for the
determination of exponents in the plane.

To be more precise, consider critical percolation with fixed mesh equal to 
$1$,
and consider the event $G_j(r,R)$ that there exists $j$
disjoint blue {\it crossings} of the semi-annulus
$A_+(r,R) := \{z \ : \  r < |z| < R,\  \Im (z) > 0  \}$.
By a blue crossing we mean a (discrete) simple blue curve,
i.e. a sequence of distinct blue hexagons
connecting the ``semicircle'' boundaries of $A_+(r,R)$,
where consecutive hexagons are adjacent.

Actually, we slightly modify the definition of $A_+ (r,R)$
as follows:
All hexagons intersecting the circle of radius $r$ will 
be in $A_+ (r,R)$ while those which intersect 
the circle of radius $R$ are not in $A_+(r,R)$.
A crossing of $A_+ (r,R)$ is then a bluse simple curve from
$|z|=r$ to a hexagon on the ``outer boundary'' of $A_+ (r,R)$.

One could as well take a semi-hexagonal or a triangular
shape instead of the semi-circles to simplify the
discrete approximation.
Define $a_j (r, R) := \Prob [ G_j(r,R)]$.

\begin{rem}[Changing the colors]\label{colorschange}
The convention to study blue crossings is not restrictive:
it is standard (see e.g., \cite {ADA,LSW1})
that for any given sequence of colours, the
probability that there exist left-right crossings
of a topological rectangle (for discrete critical percolation)
of some prescribed colours in prescribed order
is in fact independent of this sequence of
colours.
It follows that this probability is comparable (for long 
rectangles) to
the probability of $j$ crossings {\em of arbitrary colours}.
One can also check that the probabilities of {\em at least} $j$ crossings,
and of {\em exactly} $j$ crossings are comparable.
The reason is that one can explore the crossings from
``below'' one-by-one and flip all colours above some of them
without changing the probabilities of configurations.
\end{rem}

We prove the following result, which was predicted by physicists:

\begin {theorem}[Half-plane exponents]
\label {halfplane}
For any $j \ge 1$, and for all large enough $r$
(i.e. $r>\hbox{const}(j)$),
$$
a_j(r,R) = R^{- j ( j+1)/6 + o(1) }
\hbox { when } R \to \infty.
$$
\end {theorem}

We first show that the theorem is a direct consequence
of the following two facts,
which will be discussed in the remainder of the section:
\begin{itemize}
\item Identification between $SLE_6$ and continuum percolation
implies that exponents for continuum percolation
are equal to the exponents for $SLE_6$,
computed in \cite{LSW1}.
This can be written in terms of discrete percolation:
\begin{equation}\label{halfplaneexp}
\lim_{\rho\to\infty} a_j(\rho,R \rho)=
R^{- j ( j+1)/6 + o(1) }\hbox { when } R \to \infty.
\end{equation}
\item
Crossing probabilities
enjoy the following (approximate) multiplicativity property
with some positive $c=\hbox{const}(j)$,
provided $r''\ge r'\ge r>j$
(cf. \cite{Kpaper,KSZ}): 
\begin {equation}
\label {submult}
a_j(r, r') \, a_j(r' , r'') \ge a_j (r, r'') \ge c \, a_j(r, r') \, a_j(r',
r'').
\end{equation}
\end{itemize}

In fact, as $a_j (r,R)$ is decreasing in $R$,
to establish Theorem~\ref{halfplane}
it is sufficient to show that for any fixed $\epsilon$, there
exists $K>1$ such that for sufficiently large $n$,
\begin {equation}
\label {sufficient}
(K^n)^{- j (j+1)/ 6 - \epsilon} \le
a_j (r, r K^n) \le (K^n)^{-j (j+1)/6 + \epsilon} .
\end {equation}
To prove (\ref{sufficient}), we
use (\ref{halfplaneexp})
and choose large enough $K$ so that
$$
c^{-1} K^{- j (j+1)/6 - \epsilon/ 2 } \le a_j(\rho,\rho K) \le
K^{-j (j+1)/ 6 + \epsilon/ 2},
$$
for sufficiently large $\rho$.
Together with (\ref {submult}) this implies
\begin {eqnarray*}
\hbox{const} \prod_{m=1}^n K^{- j (j+1)/6 - \epsilon/ 2 }
\le \prod_{m=1}^n c a_j (K^{m-1}r, K^m r)
\le a_j (r, rK^n)\le\\
\le \prod_{m=1}^n a_j ( K^{m-1}r, K^mr)
\le \hbox{const} \prod_{m=1}^n K^{- j (j+1)/6 + \epsilon/ 2 },
\end {eqnarray*}and (\ref {sufficient}) readily follows.
To prove Theorem~\ref{halfplane}
it remains to check (\ref{halfplaneexp})
and (\ref{submult}).

\subsection{Chordal processes}

\paragraph{Chordal exploration process.}
Suppose for a moment that $\Omega$ is a simply connected set of hexagons,
and that $a$ and $b$ are two distinct vertices (of the honeycombe lattice)
that are on its boundary.
For convenience, we will often identify an arbitrary domain $\Omega$
with its discrete hexagonal approximation.
Colour all hexagons on the boundary of $\Omega$ that
are between $a$ and $b$ in counter-clockwise order (resp.
clockwise order) in blue (resp. yellow) and
call this set of blue
hexagons $\partial_b$ (resp. $\partial_y$). There exists a unique
curve separating the blue cluster attached to $\partial_b$ from
the yellow cluster attached to $\partial_y$.
This is the exploration process from $a$ to $b$ in $\Omega$.
Note that this curve is a simple curve that has blue hexagons to
its left and yellow hexagons to its right (if seen from $a$ to $b$),
and that it can be defined dynamically, as an ``exploration
process'' that turns right when it meets a yellow hexagon and
left when it meets a blue hexagon. In particular, for a fixed
curve $\gamma$, the event that the exploration process is equal
to $\gamma$ depends only on the state of the hexagons that
are in the neighbourhood of $\gamma$.

By Remark~\ref{colorschange}, $a_j(r,R)$ is also
the probability of the event of $j$ crossings of alternate colors,
which is can be described as
discrete exploration process inside the semiannulus $A_+(r,R)$
from the point $r$ (i.e., the real point $r$ in the 
complex plane) to the point $-R$ making $j$ crossings
before hitting the interval $[-R,-r]$.
Note that for odd $j$ this can also be described
as the existence of $(j-1)/2$ disjoint
yellow clusters that cross the semi-annulus.

\paragraph{Chordal $SLE_6$.}
One can view chordal $SLE_6$ as follows (e.g., \cite {LSW1,RS}):
For any simply connected open set
$\Omega\subset \C$, $\Omega \not= \C$
  and two distinguished points
(or ends) $a$ and $b$ on its boundary, it is a
  random continuous curve $\gamma$
from $a$ to $b$ in $\overline \Omega$.
The law of this curve is conformally invariant by construction,
this curve has double-points but no ``self-crossings,''
and is of fractal dimension $7/4$ almost surely (\cite {RS,B}).

Critical exponents
associated to the $SLE_6$ curve in a semi-annulus have been
computed in 
\cite {LSW1,LSW3}:
Define the probability $a_j^{sle}(r,R)$
that $SLE_6$ from $r$ to $-R$ in the
semi-annulus $A_+(r,R)$ makes $j$ crossings before its hitting time 
$\tau$ of the  interval
$[-R,-r]$.
Note also that because of the conformal invariance,
$a_j^{sle}(r,R)$ depends on the ratio $R/r$ only:
$a_j^{sle}(r,R)=a_j^{sle}(R/r)$.
 Let $\sigma$ denote the time at which $SLE_6$
performs its first crossing of the semi-annulus. 
At this time, the $SLE_6$ has still to perform $j-1$
crossings between the two circles
 in the connected component $U$ of $A_+ (r, R) \setminus 
\gamma [0,\sigma]$ before $\tau$. 
By conformal invariance and the strong Markov property of $SLE_6$ at time
$\sigma$, it follows that given $\gamma [0,\sigma]$, the conditional probability that $SLE_6$  
makes the $j$ crossings is 
$a_{j-1}^{sle} (S)$ where $ S / \pi $ 
is the extremal distance between the two circles in $U$.
Hence, 
$$    
a_j^{sle} (R/r) = 
\expect  [1_{\sigma < \tau}  a_{j-1} (S) ]
.$$
Theorem 2.6 in \cite {LSW1} states that  
that for all non-negative  $\lambda$,
$$
\expect [1_{\sigma<\tau } S^{- \lambda} ] 
\asymp
R^{- u( \lambda)  } 
~\mbox{when}~R \to \infty,
$$
where $ u ( \lambda ) = ( 6 \lambda + 1 + \sqrt {24 \lambda +1}) / 6 $,
and $\asymp$ means that the ratio
between both quantities remain bounded
and bounded away from zero.
As $u^{\circ j} (0) = j (j+1) / 6$, 
it therefore follows by induction that for all $j \ge 1$,
\begin{equation}
\label{sleexp}
a_j^{sle} (R) \asymp R^{-j(j+1)/6 },
~\mbox{when}~R \to \infty.
\end{equation}

\paragraph{Exploration process and $SLE_6$.}
In \cite {S}
it is shown among other things that
the discrete exploration process from $a$ to $b$
in the discretized approximation of $\Omega$ converges
in law towards $SLE_6$, with respect to the Hausdorff
topology on simple curves,
when the mesh of the lattice goes to zero.
In fact \cite {Sprep}, it is also possible to derive a
slightly stronger statement that loosely speaking
the convergence takes place uniformly with respect to
the domain $\Omega$ and the location of the endpoints $a$ and $b$.

Hence it follows from \cite {S}, that when the mesh of the lattice
goes to zero, the probability that the discrete exploration
process makes $j$ crossings of a semiannulus converges to that for $SLE_6$.
Alternatively we can
increase the domain while preserving the mesh, and conclude that
\begin {equation}
\label {limit}
\lim_{\rho\to\infty}a_j(\rho r, \rho R)=a_j^{sle}(r,R).
\end{equation}
Let us stress  that this uses the following
``a priori bounds''  (cf. \cite {LSWprep,Kpaper,KSZ}):
For some fixed $r$ and $\epsilon>0$ and for all sufficiently large $R$,
\begin {equation}
\label {3arms}
a_3 (r,R) \le R^{-1- \epsilon},
\end {equation}
as well as an analogous result
for 6 arms in the plane 
(see (\ref {6arms}) in the next section).
These bounds
for instance prevent the possibility that with positive
probability, a crossing of the rectangle
appears in the
scaling limit while it was not present in the discrete case,
or that two distinct crossings collapse into one
in the scaling limit
(see \cite {Sprep} for more details).

Equations (\ref{sleexp}) and (\ref{limit})
clearly imply the desired (\ref{halfplaneexp}).

\subsection{Multiplicativity}

The left inequality in the ``approximate multiplicativity''
(\ref{submult}) is immediate since $A_+ (r,r')$ and $A_+ (r', r'')$ are
disjoint.

The right inequality is elementary using an argument involving the
Russo-Seymour-Welsh theory (referred to as RSW theory in
the sequel).
If  $2 r<r'<r''/2$, one can use
the Harris-FKG (Fortuin-Kasteleyn-Ginibre)
inequality and the events that there exist j disjoint blue
crossings
of $A_+(r'/ 2, 2 r' )$, and also $j$
disjoint blue crossings of the ``long side'' of the semi-annuli
$A_+ (r' / 2 , r')$ (i.e. that join the two real segments) and
$A_+ (r', r' /2)$).
This argument relies on the fact that all crossings have the same
colour. As we shall see later, things are more involved for plane exponents
because of the fact that we will be considering crossings of the two
colours.

This is the only case needed in the proof, but remaining cases are also
easy: if e.g. $r'\ge r''/2$ then by the reasoning above
and standard RSW theory,
$$a_j(r,r'')\ge a_j(r,2r')\ge  c \, a_j(r, r')\,a_j(r',2r')
\ge  c' \, a_j(r, r') \, a_j(r',r'').$$

\section {Plane exponents}

We now study the plane exponents, i.e.
probabilities of crossings of annuli instead of semi-annuli.
There is a profound difference with the half-plane case,
when we had  a ``starting half-line,''
which allowed to choose for a given configuration
a {\em ``canonical realization''} of crossings
(the lowest, the second from below, and so on)
and to change colours of crossings at will.

In contrast, when one studies $j$ crossings of an annulus
of the same colour, then there is no canonical way to choose
their realization.
There was a prediction by theoretical physicists
 in the case $j=1$ (this is the ``one-arm
exponent'' that we will discuss later),
but not for $j\ge 2$.
However, conformal invariance
and $SLE_6$ can be used (see  \cite {LSWprep}) to describe the 
``backbone exponent'' i.e. the case  $j=2$,
which is responsible for the dimension of the ``backbone'' --
the sites of percolation cluster connected to the boundary of the domain
by two disjoint blue curves -- 
as the leading eigenvalue of a certain differential 
operator.

It turns out (as observed by \cite {ADA})
that the probabilities and the exponents
for $j$ blue crossings are different from those
for $j$ crossings which are not all of the same colour,
and for the latter case there are physical predictions,
which we are going to establish.
The reason is that one can take two neighboring crossings
of different colours, choose their closest (to each other)
realizations, and then use their hull as a starting point,
choosing ``canonical'' realizations of other crossings
and changing their colours if needed.
Similarly to the half-plane case,
one concludes that exponents are the same
for $j$ crossings of any prescribed colours in any prescribed order,
as long as colours are not all the same.

To be more precise, consider critical percolation with fixed mesh equal to 
$1$,
and consider the event $H_j(r,R)$ that there exists $j$
disjoint crossings of the annulus
$A(r,R) := \{z \ : \  r < |z| < R \}$ (in fact we modify $A (r,R)$
as $A_+ (r,R)$ in the previous section),
not all of the same colour.
Define $b_j (r, R) := \Prob [ H_j(r,R)]$.
One can also prescribe colours of the crossings and their order,
which will change $b_j$ up to a multiplicative constant (we will 
justify this rigorously later),
preserving the theorem below.
The goal of this section is to establish
the following result, which was predicted by physicists:

\begin {theorem} [Plane exponents]
\label {plane}
For any $j \ge 2$, and for all large enough $r$
(i.e. $r>\hbox{const}(j)$),
$$
b_j(r,R) = R^{- (j^2-1)/12 + o(1) }
\hbox { when } R \to \infty.
$$
\end {theorem}
The statement above includes 
(\ref {4arms}) as a particular case:
for even $j=2k$, we can take
crossings of alternate colours, and that
corresponds to the existence 
of $k$ disjoint clusters that 
cross the annulus. In particular,
$\Prob(A^2_R)=b_{4}(2,R)=R^{-5/4+o(1)}$.

Exactly as its half-plane counterpart, the theorem follows
from the two observations which will be discussed below:
\begin{itemize}
\item Identification between $SLE_6$ and continuum percolation
implies that exponents for continuum percolation
are equal to the exponents for $SLE_6$,
computed in \cite{LSW2}.
This can be written in terms of the crossings
probabilities for the exploration process
in the annulus (to be defined below):
\begin{equation}\label{planeexp}
\lim_{\rho\to\infty} b_j^{ep}(\rho,R \rho)=
R^{- (j^2-1)/12 + o(1) }\hbox { when } R \to \infty.
\end{equation}
\item
There is an unbounded set ${\cal R} \subset (0, \infty)$ and a positive $c$
(depending on $j$ only)
such that crossing probabilities
enjoy the following
approximate multiplicativity property for any $R\in{\cal R}$
and all $n \ge 1$
\begin {equation}
\label {planemult}
\mbox{const}(R)\,c^{-n} \,\prod_{l=1}^{n}b_j^{ep}(2 R^{l}, R^{l+1})
\ge b_j(R, R^{n+1})
\ge \mbox{const}(R)\,c^n \, \prod_{l=1}^{n}b_j^{ep}(2 R^{l},
 2R^{l+1}).
\end{equation}
A more general inequality, analogous to (\ref{submult})
is valid, but for simplicity we prove the above version only,
which is sufficient to establish the theorem.
\end{itemize}

\subsection{Processes in an annulus}

\paragraph{Exploration process and $SLE_6$ in the universal cover.}

Suppose that an annulus $A=A(r,R)$ is given,
and denote by $\tilde A=\tilde A(r,R)$ its universal cover,
with inherited lattice structure.
Fix a point $x$ on the outer circle of $A$,
let $\tilde x$ be one of its lifts to $\tilde A$,
and $y$ be a
``counterclockwise point at infinity on the boundary of $\tilde A$.''

Perform chordal exploration process and $SLE_6$
from ${\tilde x}$ to $y$ in the domain $\tilde A$.
We define a disconnection time $T$,
which is the first time when the projection of
the trajectory to the annulus $A$ disconnects
the inner circle from the outer one.

As discussed, the law of the chordal exploration process
converges (as mesh goes to zero) to that of $SLE_6$.
If two trajectories of the exploration process
are $\delta$-close in the Hausdorff metric,
but have drastically different disconnection times,
then ``six arms'' must occur for one of them.
The following ``a priori bound''  (cf. \cite {Aiz,KSZ})
\begin {equation}
\label {6arms}
b_6 (r,R) \le \mbox{const}\,({r}/{R})^{-2-\epsilon},
\end {equation}
imply that they occur somewhere ``at scale $\delta$''
with probability $o(1)$, $\delta\to0$.
Therefore the law of the chordal exploration process
stopped at the disconnection time
converges (as mesh goes to zero) to that of $SLE_6$
stopped at the disconnection time.

As before, we infer that the probability $b_j^{ep} (r, R)$
that the exploration process makes $j-1$ crossings between the
inner and outer
boundaries of $\tilde A$ before time $T$ converges
to the similar probability $b_j^{sle} (r, R)$ for $SLE_6$
as the mesh of the lattice goes to zero.
Increasing the domain while preserving the mesh instead, we conclude that
\begin {equation}
\label {planelimit}
\lim_{\rho\to\infty}b_j^{ep}(\rho r, \rho R)=b_j^{sle}(r,R).
\end{equation}

\paragraph{Projecting from the universal cover to the annulus.}
The mentioned chordal processes in $\tilde A$
up to the disconnection time $T$
can be projected to the annulus $A$.
Locally their definitions coincide with
the processes in the annulus
described below
(which are well-defined up to the disconnection time).
So by the restriction property
(laws of the exploration process and $SLE_6$
depend only on the neighbourhoods of their traces)
we conclude that the projections
of chordal processes in $\tilde A$
coincide with the following processes in $A$
up to the disconnection time $T$:

\begin{itemize}
\item
Exploration process in an annulus
follows the same ``blue to the right-yellow to the left''
rule as the chordal exploration process,
except that we colour the
hexagons of the inner circle in yellow, and that when the
exploration process hits the outer circle and the
continuous determination of the argument of the
exploration process is larger (resp. smaller) than that
  the starting point of the exploration, the boundary
point that it hits on the outer circle is blue (resp. yellow).
\item
$SLE_6$ in an annulus goes as radial $SLE_6$ (see \cite{LSW2})
from $x$ to the center of the annulus up to the first
hitting $\rho$ of the inner boundary.
Afterwards it continues like chordal $SLE_6$ in
the remaining domain,
until the disconnection time.
\end{itemize}

Particularly, we conclude that probabilities
$b_j^{ep}$ and $b_j^{sle}$
are the same for chordal processes in $\tilde A$
and their ``annular'' counterparts in $A$.
Note that scaling implies that just as in the 
chordal case, 
$b_j^{sle} (r,R) =: b_j^{sle} (R/r)$.

\paragraph{Exponents for $SLE_6$.}

The computation of exponents for radial and chordal $SLE_6$ in
\cite {LSW1, LSW2}
yields that for $j \ge 2$ 
\begin{equation}
\label{planesleexp}
b_j^{sle} (R) = R^{-(j^2-1)/12  + o(1) },
~\mbox{when}~R \to \infty.
\end{equation}
Indeed, if $\sigma$ now denotes the first hitting time of the inner circle
by the $SLE_6$ in the annulus $A(1,R)$, and $S / \pi$ the extremal distance
between the two circles in $A(1,R) \setminus \gamma [0,  \sigma]$, then
conformal invariance and the strong Markov property show that
$$
b_j^{sle} (R ) = \expect [ 1_{\sigma< T} a_{j-2}^{sle} (S) ].
$$
It is shown in \cite {LSW1}, Theorem 3.1 
that for 
$\lambda \ge 1$ and $\lambda =0$ (see \cite {B} for this case), 
\begin {equation}
\label {nulambda}
\expect [ 1_{\sigma < T } S^{-\lambda} ] 
\asymp 
R^{- \nu (\lambda)}~\mbox{when}~R \to \infty,
\end {equation}
where $ \nu (\lambda ) = ( 4 \lambda + 1
+ \sqrt { 1 + 24 \lambda } ) / 8$.
(\ref {planesleexp}) for $j \not= 3$
then follows by plugging in (\ref {sleexp}).  
So far, a direct proof of (\ref {nulambda}) for $\lambda=1/3$
is  missing in the literature (the proof for $\lambda \ge 1$ in \cite {LSW2} 
uses the computation of a ``derivative exponent'' and a convexity 
argument. The latter does not work directly for $\lambda \in (0,1)$).
Equation (\ref {planesleexp}) for $j=3$ can however be 
derived via other rather convoluted means (for instance 
a universality argument and analyticity of intersection 
exponents \cite {LSWa,LSWfr}). 

\begin{rem}
The fact that $b_2^{sle}(R) = R^{-1/4 + o (1)}$ is
related to the fact that the Hausdorff dimension of the $SLE_6$ curve 
is $2- 1/4=7/4$ (see \cite {RS, B}).
The fact that
$b_3^{sle} (R) = R^{- 2/3 + o(1)}$ is related to the fact that
the dimension of the outer frontiers of $SLE_6$ and of planar Brownian motion 
is $2-2/3=4/3$ (see \cite {LSWfr})  and that the fact that
$b_4^{sle}(R) = R^{-5/4 + o(1)}$ is related to the fact that the Hausdorff
dimension of the set of (local) cut points of $SLE_6$ 
is $2- 5/4 = 3/4$ (\cite {LSW2}), e.g. using the
simple identification between Brownian
hulls and $SLE_6$ hulls \cite {W}.
\end{rem}

Combining (\ref{planelimit}) and (\ref{planesleexp})
we arrive at (\ref{planeexp}).

\begin {rem}
Before proceeding, we want to remark
that events corresponding to $b_j^{ep}$ can be easily described
in terms of percolation crossings.
Namely at the time $T$, when exploration process started from the point $a$
first hits the inner circle, its hull $K$ is bounded
by two crossings of the annulus:
the clockwise-most blue crossing
and the counterclockwise-most yellow crossing,
containing $a$ between them.
After that the exploration process continues as the chordal
process in $A\setminus K$, creating $(j-2)$ more crossings
of alternate colours.
Just as in the halfplane case we can change the colours of all crossings
except the first two, so we conclude that
$b_j^{ep}$ gives (up to a multiplicative constant)
the probability of having $(j-2)$ crossings
of some prescribed (any) colours outside the hull
of two crossings of opposite colours containing $a$ between them.

So, the description of $b_j^{ep}$ adds an additional requirement
(a prescribed starting point separating two crossings of opposite colour)
as compared to $b_j$, and we conclude that
\begin{equation}\label{bepb}
b_j^{ep}(r,R)\le b_j(r,R).
\end{equation}
Unlike the half-plane case, the reverse inequality
is valid up to a multiplicative constant only.
This is trivial when $j=2k$ is even (one just has
to take alternating colours as prescribed order and then the
starting point of the exploration process is anyway between
two crossings of different colours) which corresponds to the
probability of $k$ disjoint blue clusters. But an additional
argument, in the spirit of the discussion below,
 is needed for odd number of crossings.
\end {rem}

\subsection{Multiplicativity}

Because of the two different colours,
the simple argument based on the FKG inequality and
the RSW theory can not be immediately applied.
Nevertheless, some more elaborate approaches
suitable for similar problems were developed
by Kesten and others.
Needed results for $4$ and $5$ crossings can be found
in \cite{Kpaper} and \cite{KSZ} correspondingly,
but there seems to be no readily available reference
for an arbitrary number of arms, so
we present a proof below.
Such arguments are also very close to Lawler's 
separation Lemmas for Brownian intersection
probabilities, see e.g., \cite {LSWup}.

It follows from \cite{Sprep}, that
$b_j(\rho r,\rho R)$ has a scaling limit,
which is conformally invariant, and so depends on the ratio $R/r$ only:
\begin{equation}\label{bjlimit}
b_j'(R/r)=\lim_{\rho\to\infty}b_j(\rho r, \rho R).
\end{equation}
By standard RSW theory, $b_j(r,R)$ is bounded from below
by a power of $R/r$, hence
$$b_j'(R)\ge \mbox{const}\,R^{-\zeta},$$
for some $\zeta>0$.
Therefore we can conclude that
there exist a positive constant $K$ and
an unbounded set ${\cal R}$ of radii $R$
such that
\begin{equation}\label{b4r}
b_j'(R/8)\le K b_j'(R/2).
\end{equation}

We now define a notion of $\delta$-good configurations,
when landing points of crossing are ``well separated''.
We say that a configuration is $\delta$-good
in the annulus $A(r, r')$ (we assume that $r' \ge 4r$) if
there exist $j$ disjoint crossings
not all of the same colours,
and there is no ball of radius $\delta r$
(resp. $\delta r'$) centered on the
inner (resp. outer) circle and intersecting at least three of the crossings.
We call $b_j^{good}(r,r')=b_j^{good} (r, r' , \delta)$ the probability of 
this
event.
Note that the event corresponding to
$b_j (r, r') - b_j^{good}(r, r', \delta)$ is contained in the event
that there exist $j$ crossings of the annulus ${A} ( 2r, r'/2)$ and
that three crossings
of ${A} (r, 2r)$ (or of ${A} (r'/2, r')$)
come $\delta r$-close (resp. $\delta r'$ close) near the inner
circle (resp. outer circle) of the annulus.
From  (\ref {3arms}) it follows that the probabilities
of these last events go to zero when $\delta \to 0$,
(uniformly with respect to $r$ and $r'$)
so we can fix $\delta:=\mbox{const}(K)$,
so that these probabilities are smaller that $1/(8K)$
and hence for all $r'/r \ge 4$,
\begin{equation}\label{badgood}
b_j^{good} ( r, r' , \delta)
\ge b_j(r, r' ) - \frac1{4K}\,b_j(2r, r'/2).
\end{equation}

Fix $R\in{\cal R}$.
It follows from (\ref{bjlimit}) and (\ref{b4r})
that there exist $l_R:=\mbox{const}(R)$ such that for
$l\ge l_R$
\begin{equation}\label{double}
b_j\left(4R^{l},R^{l+1}/2\right)\le 2K\,b_j\left(2R^{l},R^{l+1}\right).
\end{equation}
Combining (\ref{badgood}) and (\ref{double}) we conclude that for
$l\ge l_R$
\begin{equation}\label{good}
b_j^{good} (2R^{l}, R^{l+1})
\ge\left(1-\frac1{4K}\,2K\right)\, b_j(2R^{l}, R^{l+1} )
=\frac12\,b_j(2R^{l}, R^{l+1}).
\end{equation}

Standard (but delicate) 
techniques based on RSW theory (cf. \cite{Kpaper})
show that there exists a constant $Q=\mbox{const}(\delta) \ge 1$
(note that $Q$ depends on $\delta=\mbox{const}(K)=\mbox{const}(j)$),
such that
\begin{equation}\label{goodmult}
b_j^{good}(r,r')\,b_j^{good}(2r',r'')\le Q\, b_j^{good}(r,r'').
\end{equation}
Similarly, the exists a constant $Q'\ge 1$ such that 
\begin {equation}\label{goodep}
b_j^{good} (r, r'/2) \le Q' \, b_j^{ep} (r,r').
\end {equation}
We are now ready to conclude,
writing for the left half of (\ref{planemult})
$$b_j (R, R^{n+1}) 
\le \prod_{l=l_K}^n b_j (2R^l, R^{l+1}) 
\stackrel {(\ref {good})}{\le}
 2^n \prod_{l=l_K}^n b_j^{good} (2R^l, R^{l+1}) 
\stackrel {(\ref {goodep})}{\le}
 C(R)(2Q')^n \prod_{l=1}^n b_j^{ep} (2R^l, 2R^{l+1})
$$
(in the first inequality we used
that the event corresponding to $b_j(R, R^{n+1})$
requires simultaneous occurrence of independent events
corresponding to $b_j(2 R^{l}, R^{l+1})$ with $l=l_K,\ldots, n$).
For the right half of (\ref{planemult}) we write
\begin {eqnarray*}
b_j(R, R^{n+1}) \ge  
b_j^{good}(R, R^{n+1})
\stackrel{(\ref{goodmult})}{\ge}
C(R)
Q^{-n}\prod_{l=l_K}^{n}b_j^{good}(2 R^{l}, R^{l+1}) \ge \\
\stackrel{(\ref{good},\ref{bepb})}{\ge}
C(R)
(2Q)^{-n}\prod_{l=1}^{n}b_j^{ep}(2 R^{l}, R^{l+1}).
\end {eqnarray*}

\subsection {One crossing of the annulus}

In order to derive
(\ref {onearm}), one has
to translate the existence of one blue connection between circles
in terms of the exploration processes.
Consider as in the 
previous subsections
 a discrete exploration process in the annulus $A(r,R)$ up to its
disconnection time $T$ i.e. the first time at which the exploration process
contains a closed loop around the inner circle. Let $\rho$ denote the 
first hitting time of the inner circle.
If $\rho<T$, then it means that there exists one arm of each colour 
joining the two circles, and in particular a blue one. If $T< \rho$,
then one has to see in which direction $\gamma$ did wind around 
the inner circle: If $\gamma [0,T]$ contains a clockwise loop
around the inner circle, then it means that the 
exploration process has discovered a closed loop of yellow hexagons
around the inner circle, and in this case, there is no blue
connection between the inner and the outer circle.
If however $\gamma [0,T]$ makes an anti-clockwise loop
around the inner circle, then the exploration process has 
discovered a closed loop $l$ of blue hexagons around the 
inner circle that is connected to the outer circle by a blue
path. Furthermore, the exploration process has not explored any
of the hexagons that are in the connected component of 
$A(r,R) \setminus l$ containing the inner circle.
Hence, to see if there is a blue crossing of the annulus, it 
remains to see if there is a blue crossing between the inner
circle and $l$ i.e. to start the same algorithm again in this
new domain.

Hence, one is lead to study the following
quantities: What is the probability that the 
radial $SLE_6$ in the unit disc up to its first hitting
time $\rho$ of the circle of radius $r$ contains no clockwise 
loop around the origin? 
In \cite {LSWprep}, this probability is shown to decay like
$r^{5/48}$.
The number $5/48$ in fact corresponds to the same differential operator
as that describing the probability of no loop at all (i.e. corresponding
to the exponent $1/4$), but with different
 boundary conditions (one Dirichlet and one Neumann instead of two Dirichlet). 

Using arguments in the same spirit than those described above, 
one can then show that indeed,  
$\Prob [ A_R^1 ] = R^{-5/48 + o (1) }$.
For details, we refer to \cite {LSWprep}.
Note that this is equivalent to  
$$\Prob [ 0 \hbox { is connected to } x  ] 
  = |x|^{- 5/24 + o(1)}   
$$
when $x$ goes to infinity (i.e., the exponent often denoted by $\eta$ 
exists and is equal to $5/24$).

\section {Some open questions}

To conclude, we very briefly list some questions that seem still open
at this moment. Some of them are probably within reach, and some
are less accessible.

\begin {enumerate}
\item
Generalizing the results of Smirnov \cite {S}, to other lattices. The first
two natural candidates are bond percolation on the square lattice and site
percolation
on a Voronoi tesselation (see e.g. \cite {BS}),
 that both have a ``self-duality'' type
property (in particular, the value $p=1/2$ has to be studied).
It would be sufficient to prove Cardy's formula,
but in both cases the method used in \cite {S} does not apply directly.

\item
Existence of two, three, and four arms from the vicinity of
a site represent it belonging to frontier of a percolation cluster,
perimeter of a percolation cluster, and being a
pivotal site respectively.
Thus we infer that on a lattice with mesh $1$, when we speak
of clusters of size $\approx N$,
a site has probablity $\approx N^{-1/4}$
to belong to a frontier,
probability $\approx N^{-2/3}$ to belong to perimeter,
and probability $\approx N^{-5/4}$
of being pivotal.
One should be able to show a stronger statement, 
roughly speaking that a cluster of size $\approx N$
has frontier of $\approx N^{2-1/4}=N^{7/4}$ sites,
perimeter of $\approx N^{2-2/3}=N^{4/3}$ sites,
and $\approx N^{2-5/4}=N^{3/4}$ pivotal sites.
The counterpart of this stronger claim for continuum percolation
follows from the identification
between $SLE_6$ and the scaling limit of percolation cluster perimeter:
e.g., the scaling limit of the frontier
of a percolation cluster has the same law (when properly normalized)
as the Brownian frontier,
and hence has Hausdorff dimension $4/3$ almost surely. See \cite {B}
for a more direct proof. 

\item
Show that the power laws hold up to constants. For instance, does
$\Prob [ A_R^1 ] \in [c R^{-5/48}, C R^{-5/48}]$
hold for some constants
$c,C\in(0,\infty)$?
Estimates up to constants can be useful in order to derive 
results on ``discrete fractal dimension'' (see the previous 
question).

\item
Show that the exponent $\alpha$ associated to the mean number of 
clusters per vertex exists and determine its value
(conjectured to be $-2/3$).

\item
Determine the 
exponents corresponding to $j\ge 3$ blue crossings of an annulus.

\item
Understand the relation between other critical
lattice models such as the critical random cluster models and the
relation to other $SLE_\kappa$ and their critical exponents.

\end {enumerate}

{\bf Acknowledgements.} We thank Raphael Cerf and Oded Schramm for 
valuable comments and Harry Kesten for very useful
advice.

\null

S. Smirnov, 
Department of Mathematics, 
Royal Institute of Technology (KTH), 
S-10044 Stockholm, Sweden. e-mail: stas@math.kth.se

W. Werner, D\'epartement de Math\'ematiques, B\^at. 425, 
Universit\'e Paris-Sud, F-91405 Orsay cedex, France.
e-mail: wendelin.werner@math.u-psud.fr

\end {document}